\documentclass[reqno,centertags, 12pt]{amsart}
\usepackage{pdfsync}
\usepackage{graphicx}
\numberwithin{equation}{section}

\newcommand{\vp}{\varphi}
\newcommand{\T}{\partial\mathbb{D}}

\newcommand{\ds}{\displaystyle}
\newcommand{\ol}{\overline}
\newcommand{\supp}{\text{\rm{supp}}}

\newcommand{\be}{\begin{equation}}
\newcommand{\ee}{\end{equation}}
\newcommand{\ba}{\begin{array}}
\newcommand{\ea}{\end{array}}

\newcommand{\bpm}{\begin{pmatrix}}
\newcommand{\epm}{\end{pmatrix}}

\newcommand{\mR}{\mathbb{R}}

\newcommand{\ye}{\left( \ds \frac y 2 \right)}

\newtheorem{lemma}{Lemma}[section]
\newtheorem{theorem}{Theorem}[section]

\begin{document}

\title[]
{Point mass insertion on the real line and non-exponential perturbation of the recursion coefficients}
\author[M.-W. L. Wong]{Manwah Lilian Wong}
\thanks{$^*$ Mathematics 253-37, California Institute of Technology, Pasadena, CA 91125.
E-mail: wongmw@caltech.edu}
\date{June 3, 2009}

\maketitle

\begin{abstract}
We present the construction of a probability measure $d\gamma$ with compact support on $\mathbb{R}$ such that adding a discrete pure point results in changes in the recursion coefficients without exponential decay.

\end{abstract}

\section{Introduction}

Suppose $d\mu$ is a probability measure on the unit circle $\T$. We define an inner product and a norm on $L^2(\T, d\mu)$ respectively as follows:
\begin{eqnarray}
\left\langle f ,g \right\rangle & = & \ds \int_{\T} \ol{f(e^{i \theta})} g(e^{i \theta}) d\mu(\theta) \label{innerproduct}\\
\|f\|_{d\mu} & = & \left( \ds \int_{\T} |f(e^{i \theta})|^2 d\mu(\theta) \right)^{1/2}
\label{norm}
\end{eqnarray}

Using the inner product defined above, we can orthogonalize $1, z, z^2, \dots$ to obtain the family of monic orthogonal polynomials associated with the measure $d\mu$, namely, $(\Phi_n(z, d\mu))_{n}$. We denote the normalized family as $(\vp_n(z, d\mu))_{n}$.

Closely related to $\Phi_n(z)$ is the family of reversed polynomials, defined as $\Phi_n^*(z)=z^n \ol{\Phi_n(1/\ol{z})}$. They obey the well-known Szeg\H o recursion relation
\be
\Phi_{n+1}(z) = z \Phi_n (z)- \ol{\alpha_n} \Phi_n^*(z)
\label{eq01}
\ee and $\alpha_n$ is known as the $n$-th Verblunsky coefficient. The Szeg\H o recursion relations for the normalized family of orthogonal polynomials is
\begin{align}
\vp_{n+1}(z) & = (1-|\alpha_n|^2)^{-1/2}(z \vp_n(z)-\ol{\alpha_n} \vp_n^*(z)) \label{normrec1}
%\vp_{n+1}^*(z) & = (1-|\alpha_n|^2)^{-1/2}(\vp_n^*(z) - \alpha_n z \vp_n(z)) \label{normrec2}
\end{align}

These recursion relations will be useful later in this paper.

Now we turn to a probability measure $d\gamma$ on $\mathbb{R}$. We can define an inner product and norm on $L^2(\mathbb{R}, d\gamma)$ as in (\ref{innerproduct}) and (\ref{norm}), except that in this case it does not involve any conjugation. By the Gram--Schmidt process, we can orthogonalize $1,x,x^2, \dots$ and form the family of monic orthogonal polynomials, $(P_n(x))_{n=0}^{\infty}$. Upon normalization, we obtain the family of orthonormal polynomials, $(p_n(x))_{n=0}^{\infty}$. These polynomials satisfy the following three-term recursion relation
\be
%x P_n(x)& = &P_{n+1}(x) + b_{n+1} P_n(x) + a_n^2 P_{n-1}(x) \\
xp_n(x) = a_{n+1} p_{n+1}(x) + b_{n+1} p_n(x) + a_n p_{n-1}(x)
\ee where $a_n$ and $b_n$ are real numbers with $a_n > 0$. They are called the recursion coefficients of $d\gamma$.

The main result of this paper is as follows:
\begin{theorem} There exists a purely absolutely continuous measure $d\gamma_0$ supported on $[-2,2]$ with no eigenvalues outside of $[-2,2]$, such that if we add a pure point $x_0 \in \mathbb{R}\backslash[-2,2]$ in the following manner
\be
d\tilde{\gamma}(x) = (1-\beta) d\gamma_0(x) + \beta \delta_{x_0} \quad  \quad \beta>0
\label{addx0}
\ee it will result in non-exponential perturbation of the recursion coefficients $a_n(d\gamma_0)$ and $b_n(d\gamma_0)$.
\end{theorem}

This example is of particular interest because of the following history: back in 1946, Borg \cite{borg} proved a well-known result concerning the Sturm--Liouville problem that in general, a single spectrum is insufficient to determine the potential. Later, Gel'fand--Levitan \cite{gelfand} showed that in order to recover the potential one also needs the norming constants.

Norming constants correspond to the weights of pure points and it is known that in the short range case (in orthogonal polynomials language, $a_n -1, b_n \to 0$ fast), varying the norming constants will result in exponential change in the potential.

Moreover, when considering the effect of varying the weight of discrete point masses on orthogonal polynomials (both on $\mathbb{R}$ and $\T$), Simon proved that it will result in exponential perturbation of the recursion coefficients (see Corollary 24.4 and Corollary 24.3 of \cite{simon3}).

All the results mentioned above lend to a few the intuition that if the recursion coefficients $a_n \to 1$ and $b_n \to 0$ fast, then adding a pure point will result in exponential change in the recursion coefficients. However, it turned out not to be the case!

\section{Tools Involved in the Proof}

\subsection{The Szeg\H o Mapping} \label{szego} It turns out that one can relate measures supported on $[-2,2]$ with a certain class of measures on $\T$. 

Note that the map $\theta \mapsto 2 \cos \theta$ is a two-one map from $\T$ to $[-2,2]$. Therefore, given a non-trivial probability measure $d\xi$ on $\T$ that is invariant under $\theta \rightarrow -\theta$, we can define a measure  \be
d\gamma = {\rm Sz}(d\xi)
\ee using what is known as the Szeg\H o map, such that for $g$ measurable on $[-2,2]$,
\be
\ds \int g(2\cos \theta) d\xi(\theta) = \ds \int g(x) d\gamma(x)
\ee

Conversely, if we have a probability measure $\beta$ supported on $[-2,2]$, we can obtain a probability measure 
\be
\nu={\rm Sz^{-1}}(d\gamma)
\ee on $\T$ by what is known as the Inverse Szeg\H o Mapping, such that for $h(z)$ measurable on $\T$,
\be
\ds \int h(\theta) d\nu(\theta) = \int h\left(\cos^{-1}\frac x 2 \right) d\gamma(x)
\ee

There are many interesting results about the Szeg\H o mapping (see Chapter 13 of \cite{simon2}), but the only relevant one for this paper is the following by Geronimus \cite{directgeronimus} (see also Theorem 13.1.7 of \cite{simon2}):
\begin{theorem}[Geronimus \cite{geronimus1}] Let $d\xi$ be a probability measure on $\T$ which is 
invariant under $\theta \rightarrow -\theta$ and let $d\gamma={\rm Sz}(\xi)$. 
Let $\alpha_n \equiv \alpha_n(d\xi)$, $a_n \equiv a_n(d\gamma)$ and $b_n \equiv b_n(d\gamma)$. Then for $n=0,1,2, \dots$,
\begin{eqnarray}
a_{n+1}^2 &=&(1-\alpha_{2n-1})(1-\alpha_{2n})^2(1+\alpha_{2n+1}) \label{an}\\
b_{n+1} & = & (1-\alpha_{2n-1}) \alpha_{2n} - (1+\alpha_{2n-1})\alpha_{2n-2} \label{bn}
\end{eqnarray} with the convention that $\alpha_{-1}=-1$.
\label{directgeronimus}
\end{theorem}

\subsection{The Point Mass Formuula}
We add a point mass $\zeta = e^{i\omega} \in \T$ with weight $0<\beta<1$ to $d\mu$ in the following manner:
\be
d\nu = (1-\beta) d\mu + \beta \delta_\omega
\label{dnudef}
\ee  Our goal is to investigate $\alpha_n(d\nu)$. 

Point mass perturbation has a long history (see the Introduction of \cite{wong1}). One of the classic results is the following theorem:
\begin{theorem}[Geronimus \cite{geronimus1, geronimus2}] Suppose the probability measure is defined as in (\ref{dnudef}). Then
\be \Phi_{n}(z, d\nu) = \Phi_n(z) - \ds \frac{\vp_n(\zeta) K_{n-1}(z, \zeta)}{(1-\beta) \beta^{-1} + K_{n-1}(\zeta, \zeta)}
 \label{geronimus}
\ee where \be K_{n}(z,\zeta)  = \ds \sum_{j=0}^{n} \ol{\vp_j(\zeta)} \vp_j(z) \label{kndef} \ee and all objects without the label $(d\nu)$ are associated with the measure $d\mu$. 
\end{theorem}

Since $\Phi_n(0) = -\ol{\alpha_{n-1}}$, by putting $z=0$ into (\ref{geronimus}) one gets a formula relating the Verblunsky coefficients of $d\mu$ and $d\nu$. For more on the formula (\ref{geronimus}), the reader may refer to Nevai \cite{nevai, nevai1}, and Cachafeiro--Marcell\'an \cite{cm1, cm2, cm3, cm4, cm5}.

Using a totally different approach, Simon \cite{simon2} found the following formula for OPUC:
\be
\alpha_{n}(d\nu)=\alpha_{n}-q_{n}^{-1} \beta \ol{\vp_{n+1}(\zeta)} \left( \ds \sum_{j=0}^{n} \alpha_{j-1} \frac{\|\Phi_{n+1}\|}{\|\Phi_j\|}\vp_j(\zeta) \right)
\label{simonformula}
\ee where $q_{n}  =  (1-\beta) + \beta K_{n}(\zeta) ; \alpha_{-1}  = -1$.

Simon's result lays the foundation for the point mass formula. In \cite{wong2, wong1}, Wong applied the Christoffel--Darboux formula to (\ref{simonformula}) and within a few steps from (\ref{simonformula}) proved the following formula for $\alpha_n(d\nu)$:
\be
\alpha_n(d\nu) = \alpha_n(d\mu) + \Delta_n(\zeta)
\label{addpoint}
\ee where
\be
\Delta_n(\zeta) = \ds \frac{(1-|\alpha_n|^2)^{1/2} \ol{\vp_{n+1}(\zeta)}\vp_n^*(\zeta)}{(1-\beta)\beta^{-1} + K_n(\zeta)}; \quad K_n(\zeta) = \ds \sum_{j=0}^{n} |\vp_j(\zeta)|^2
\label{deltandef}
\ee 

Formula (\ref{deltandef}) turns out to be very useful (see for example, \cite{wong2, wong1, wong3}).

\section{Outline of the Proof}

\subsection{Case 1: $x_0>2$} We construct a measure $d \gamma_0$ with recursion coefficients $(a_n)$ and $(b_n)$ satisfying 
\begin{eqnarray}
a_n \nearrow 1  \quad \quad b_n \equiv 0 \label{eqn1a}\\
\ds \sum_{n} |a_n-1|^2 = \infty \label{eqn2a}
\end{eqnarray}

The measure $d\gamma_0$ is purely absolutely continuous and symmetrically supported on $[-2,2]$, with no pure points outside $[-2,2]$. We scale it by a factor $0<y<2$ to form the measure $d\gamma_y$ supported on $[-y,y]\subset [-2,2]$ (we will show the connection between $y$ and $x_0$ a bit later; see (\ref{ychoice})).

Then we use the Inverse Szeg\H o map on $d\gamma_y$ to obtain the measure $d\mu_y$. By looking at the Direct Geronimus Relations (\ref{an}) and (\ref{bn}), we find necessary conditions for $\alpha_n(d\mu_y)$ so that both (\ref{eqn1a}) and (\ref{eqn2a}) hold.

Since $d\gamma_y$ is supported on $[-y,y] \subset [-2,2]$, we know that $d\mu_y$ is supported on two identical bands. Besides, $d\mu_y$ is symmetric along both the $x$-and $y$-axes because of the symmetry of $d\gamma_y$ and the Szeg\H o map.

We add a pure point at $z=1$ to $d\mu_y$ to form the measure $d\tilde{\mu}_y$ and compute the perturbed Verblunsky coefficients $\alpha_n(d\tilde{\mu}_y)$.

Then we use the Szeg\H o map on $d\tilde{\mu}_y$ to obtain the probability measure $d\tilde{\gamma}_y$ on $\mathbb{R}$. Finally, we scale $d\tilde{\gamma}_y$ to form the measure $d\tilde{\gamma}$.

Note that if we have chosen $y$ such that
\be
\ds \frac y 2 = \ds \frac{2}{|x_0|}
\label{ychoice}
\ee then we have $d\tilde{\gamma} = (1-\beta) d\gamma_0 + \beta \delta_{x_0}$.

As the final step, we show that for some constants $C_{x_0}, D_{x_0}$ (both dependent on $x_0$) such that
\begin{eqnarray}
a_n(d\tilde \gamma) & = & a_n(d\gamma_0) + \ds \frac {C_{x_0}} {n^{3/2}} + o \left( \ds \frac 1 {n^{3/2}} \right) \\
b_n(d\tilde{\gamma}) & = & b_n(d\gamma_0) +  \ds \frac {D_{x_0}} {n^{3/2}}+ o \left( \ds \frac 1 {n^{3/2}} \right)
\end{eqnarray}

\subsection{Case 2: $x_0 < -2$} Everything in Case 1 will follow except that we add a point $z=-1$ to $d\mu_y$ instead. As we shall see later in the proof, $d\mu_y$ is symmetric both along the $x-$ and $y-$ axes. Therefore, adding a pure point at $z=-1$ is the same as adding a pure point at $z=1$ and then rotating the measure by an angle of $\pi$.

\vspace{2cm}

For the convenience of the reader, here is a diagram showing all the measures involved. We will start from the measure $d\mu_y$, and move along two directions:
\be
d\gamma_0 \overset{\text{scaling}}{\longleftarrow}
 d\gamma_y \overset{\text{Sz}^{-1}}{\longleftarrow} 
 d\mu_y \overset{\text{add} \,\, z=1}{\longrightarrow} d\tilde \mu_y  \overset{\text{Sz}^{-1}}{\longrightarrow} d\tilde\gamma_y\overset{\text{scaling}}{\longrightarrow} d\tilde \gamma
\ee

\section{The Proof}

Let $d\gamma_0$ be a probability measure  on $\mR$ with recursion coefficients satisfying (\ref{eqn1a}) and (\ref{eqn2a}).

This measure, supported on $[-2, 2]$, is purely absolutely continuous, and has no eigenvalues outside $[-2, 2]$. Moreover, if we write $d\gamma_0(x)=f(x) dx$, $f(x)$ is symmetric.

Now we scale $d\gamma_0$ to form the measure $d\gamma_y$ defined by
\be
d\gamma_y(x) = d\gamma\left( 2x y^{-1} \right) \quad 0 < y < 2
\ee

The measure $d\gamma_y$, supported on $[-y, y] \subset [-2,2]$, is purely absolutely continuous and the a.c. part of $d\gamma_y(x)$ is
\be
f_y(x) = f(2x y^{-1}) \chi_{[-y,y]}
\ee which is also symmetric.

It is well known that scaling has the following effects on the recursion coefficients
\be
a_n(d\gamma_y) = \ds \ye a_n(d\gamma_0) \quad b_n(d\gamma_y)= \ds \ye b_n(d\gamma_0)
\label{scaling}
\ee

Now we apply the inverse Szeg\H o map to $d\gamma_y$ to form the probability measure $d\mu_y$ on $\T$, see figure below:

%\begin{center}
%\includegraphics[scale=0.5]{graph.pdf}
%\end{center}
%\caption[testing]

\begin{figure}[h]
\centering
\includegraphics[scale=0.3]{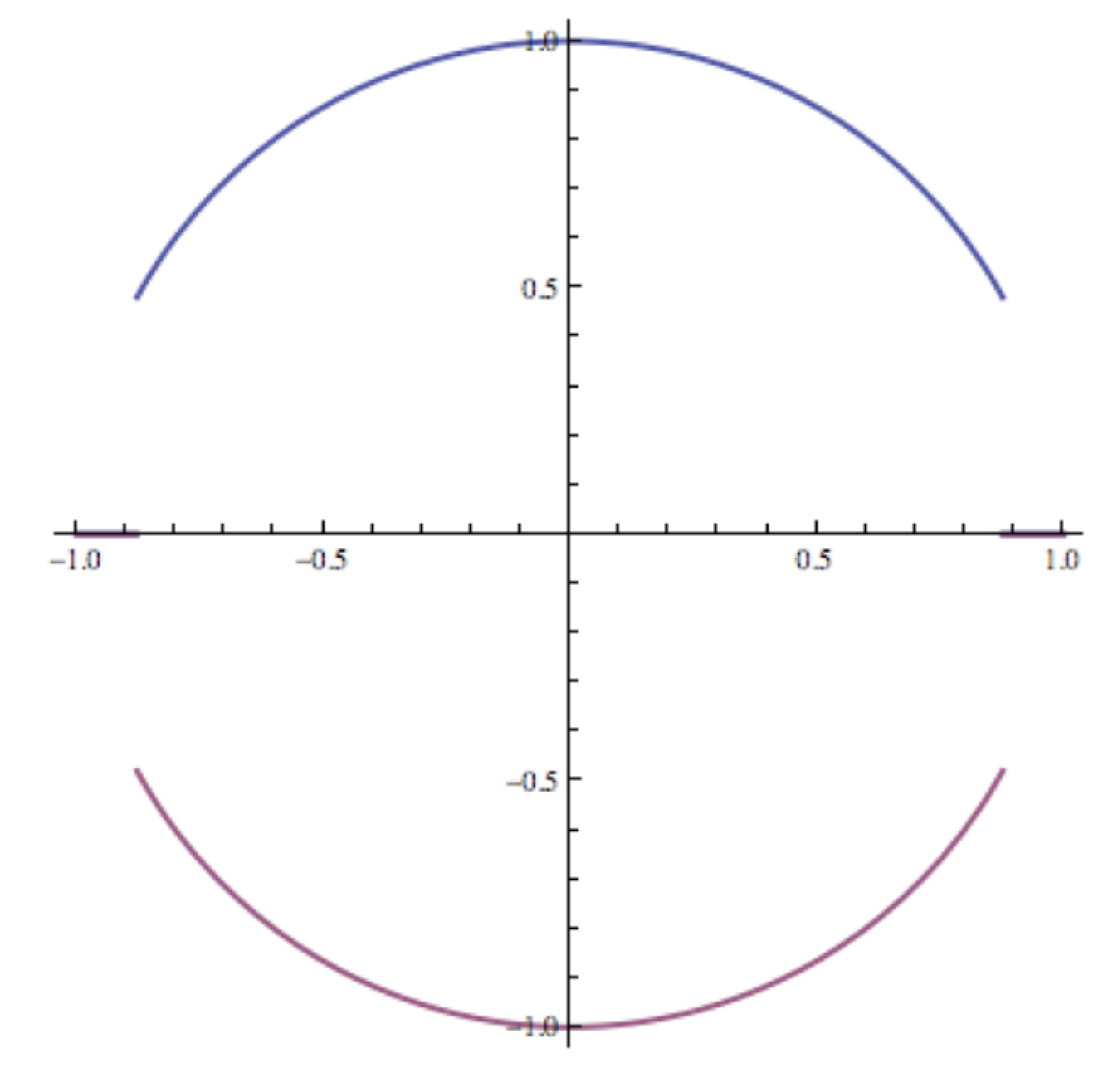}
  \caption{Graph of $\supp(d\mu)$}
\end{figure}

The measure $d\mu_y$ is supported on two arcs, $[\theta_y, \pi-\theta_y]$ and $[\pi+\theta_y, 2\pi - \theta_y]$, with a.c. part
\be
w_y(\theta) = w_y(\theta)|_{[\theta_y, \pi-\theta_y]} + w_y(\theta)|_{[\pi+\theta_y, 2\pi - \theta_y]}
\label{eqn101}
\ee where
\begin{eqnarray}
w_y(\theta) & = & 2\pi |\sin(\theta)| f_y(2\cos \theta ) \chi_{[\theta_y, \pi - \theta_y] }(\theta)
\label{eqn110} \\
\theta_y & = &\cos^{-1}\left(\frac y 2\right) \in \left(0, \frac{\pi}{2} \right)
\end{eqnarray}

By Corollary 13.1.8 of \cite{simon2}, $b_n(d\gamma_y) \equiv 0$ if and only if $\alpha_{2n}(d\mu_y) \equiv 0$. Therefore, we can express the Verblunsky coefficients of $d\mu_y$ as
\be
0, \tau_0, 0, \tau_1, 0, \tau_2, \dots
\label{eqn2}
\ee with $\tau_j = \alpha_{2j+1}$. Moreover, by Theorem 13.1.7 of \cite{simon2}, we know that
\be \begin{array}{ll}
a_{n+1}^2(d\gamma_y) & = (1- \alpha_{2n-1}(d\mu_y))(1-\alpha_{2n}(d\mu_y)^2) (1+ \alpha_{2n+1}(d\mu_y)) \\ & = (1- \tau_{n-1} ) (1+ \tau_n)
\end{array}
\label{eqn100}
\ee

Now we will choose a suitable family of $\tau_n \in \mR$ such that the corresponding $a_n(d\gamma_y)$ satisfy both (\ref{eqn1a}) and (\ref{eqn100}).

Observe that by (\ref{eqn100}) above,
\be
a_{n+1}(d\gamma_y)^2 -a_n(d\gamma_y)^2 = (1-\tau_{n-1})(\tau_n - \tau_{n-1} ) + (1+\tau_{n-1})(\tau_{n-1} - \tau_{n-2})
\ee Therefore, if we have an increasing family of $\tau_n < 0$ such that
\be
\tau_n \nearrow \tau_\infty = - \sqrt{1- \left(\ds \frac{y}{2} \right)^2}< 0
\label{e10}
\ee then $a_n(d\gamma_y) \nearrow y/2$ and the corresponding measure $d\mu_y$.

In particular, if we let
\be
\tau_k  = \ds \tau_\infty - \frac 1 {\sqrt k } %\, ;\quad y=\sqrt{3}
\ee then the goal is achieved.

Next, we prove the following lemma:

\begin{lemma} Let $d\mu_y$ be the measure on $\T$ with Verblunsky coefficients as in \eqref{eqn2} where for all large $n$,
\be
\tau_n = \tau_\infty - \ds \frac 1 {\sqrt n} \quad \quad -1<\tau_\infty<0
\ee We add a pure point at $z=1$ to $d\mu_y$ as in (\ref{dnudef}) to form $d \tilde{\mu}_y$. Then for $n=2m$ or $2m+1$, $\Delta_n(1)$ has the following expansion
\be
\Delta_n(1) = -\tau_\infty + \frac 1 {\sqrt m} + 0 + \left( 1 + \ds \frac{1}{2\tau_\infty}\right)\ds \frac 1 {m^{3/2}} + o\left( \ds \frac{1}{m^{3/2}}\right)
\ee Therefore,
\be
\alpha_n(d\mu_y) = \begin{cases}
-\tau_m + \ds \left( 1+\frac 1 {2\tau_\infty} \right) \ds \frac{1}{m^{3/2}} + e_m & n=2m \\
\ds \left( 1+\frac 1 {2\tau_\infty} \right) \ds \frac{1}{m^{3/2}} + e_m & n = 2m+1
\end{cases}
\ee where $e_m=o\left(m^{-3/2}\right)$.
\label{lemma1}
\end{lemma}

\begin{proof} Since all the Verblunsky coefficients of $d\mu_y$ are real, by induction on the recursion relation (\ref{normrec1}),
\be \begin{array}{ll}
\vp_n(1) & = \ds \prod_{j=0}^{n-1} \sqrt{\ds \frac{1-\alpha_j}{1+\alpha_j}}
\end{array}
\ee

By (\ref{eqn2}), when $n=2m$ or $2m+1$,
\be
\vp_n^*(1) = \vp_n(1) = \ds \prod_{j=0}^{m-1} \sqrt{\ds \frac{1-\tau_j}{1+\tau_j}} 
\label{vpformula}
\ee This formula will play a crucial role in the computation below.

\subsection{n is even} First, we compute $\Delta_n(1)$ when $n=2m$ using the point mass formula (\ref{deltandef}). Let
\begin{eqnarray}
A_n & =  &\ol{\vp_{n+1}(1)} \vp_{n}^*(1) \label{eqn61a} \\
B_n & =  & (1-\gamma) \gamma^{-1} + K_{n}(1,1) \label{eqn61b}
\end{eqnarray}

Then
\be
\lim_{m \to \infty} \Delta_{2m}(1) = \lim_{m \to \infty} (1-|\alpha_{2m}|^2)^{1/2} \, \ds \frac{A_{2m}}{B_{2m}} = \lim_{m \to \infty} \ds \frac{A_{2m}}{B_{2m}}
\ee because $\alpha_{2m}=0$ for all $m$. However, instead of computing this directly, we use the Stolz--Ces\`aro theorem (see \cite{cesaro}), which reads as follows
\begin{theorem}[Stolz--Ces\`aro \cite{cesaro}] Let $(\Gamma_k)_k, (\Theta_k)_k$ be two sequences of numbers such that $\Theta_n$ is positive, strictly increasing and tends to infinity. If the following limit exists,
\be
\ds \lim_{k \to \infty} \frac{\Gamma_k - \Gamma_{k-1}}{\Theta_k - \Theta_{k-1}}
\ee then it is equal to $\lim_{k \to \infty} \Gamma_k/\Theta_k$.
\label{cstheorem}
\end{theorem}

First, note that $\tau_k \to \tau_\infty < 0$. Thus,
\be
B_m \approx K_n(1,1) > |\vp_n(1)|^2 \to \infty
\ee by (\ref{vpformula}). Hence, it is legitimate for us to use Theorem \ref{cstheorem} above.

Let $K_n \equiv K_n(1,1)$ and $\vp_n \equiv \vp_n(1)$. Observe that $\vp_{2m+1} = \vp_{2m}$. Therefore, by (\ref{vpformula}),
\be
B_{2(m+1)} - B_{2m} = \vp_{2(m+1)}^2 + \vp_{2m}^2
% = \left( \ds \prod_{j=0}^{m-1} \sqrt{\frac{1-\tau_j}{1+\tau_j}}   \right) \left( \ds \frac{1-\tau_m}{1+\tau_m} + 1 \right)
 = \ds \frac{2 \vp_{2m}^2}{1 + \tau_m}
 \ee and
\be
A_{2(m+1)} - A_{2m} = \vp_{2(m+1)}^2 - \vp_{2m}^2 =  \left( \ds \frac{-2 \tau_m}{1+\tau_m}\right) \vp_{2m}^2
\ee

As a result,
\be
\lim_{m \to \infty} \Delta_{2m}(1) = \ds \lim_{m \to \infty} \frac{A_{2(m+1)}-A_{2m}}{B_{2(m+1)} - B_{2m}} = -\tau_\infty 
\ee

Next, we will prove that the rate of convergence is
\be
 \Delta_{2m}(1) = - \tau_\infty + \ds \frac{1}{\sqrt{m}} + o\left( \ds \frac{1}{\sqrt m}\right)
 \label{e0}
\ee by computing the following limit
\be
\ds \lim_{m \to \infty} m \left({\Delta_{2m}(1)+\tau_{\infty}} \right)= 1
\label{e1}
\ee

Recall the definition of $\Delta_n(1)$ and the facts $\alpha_{2m} \equiv 0$ and $\vp_{2m+1} \vp_{2m} = \vp_{2m}^2$. Thus, the left hand side of (\ref{e1}) can be expressed as $X_n / Y_n$, where
\begin{eqnarray}
X_m &=& \sqrt{m} \left[ \vp_{2m}^2 + \tau_{\infty} K_{2m} \right] \\
Y_m & =& K_{2m} \to \infty
\end{eqnarray}

We use the Stolz--Ces\`aro Theorem again. First, observe that
\be
Y_{m+1} - Y_m= \ds \frac{1-\tau_m}{1+\tau_m}+1 =\ds \frac{2}{1+\tau_m} \vp_{2m}^2
\ee

Then we compute
\begin{multline}
X_{2(m+1)} - X_{2m} \\
=\underbrace{\left[ \sqrt{m+1} \ds \frac{1-\tau_m}{1+\tau_m} - \sqrt{m} \right] \vp_{2m}^2 }_{\text{(I)}}+\underbrace{ \tau_\infty \left[ \sqrt{m+1}K_{2(m+1)} -\sqrt m K_{2m}\right]}_{\text{(II)}}
\label{e2}
\end{multline}

Consider each term in (\ref{e2}) above.
 \be  \ds \frac{\text{(I)}}{Y_{2(m+1)}-Y_{2m}} =\ds \frac{\sqrt{m+1} \frac{1-\tau_m}{1+\tau_m} - \sqrt m}{\frac{2}{1+\tau_m}}=\ds \frac{\sqrt{m+1}(1-\tau_m)-\sqrt{m}(1+\tau_m)}{2}
%& = \ds \frac 5 2 + \frac 1 m + \frac m 2 - \tau_\infty - \frac{3 m \tau_\infty} 2
\label{e4}
\ee 

Moreover,
\be
\text{(II)} = \tau_\infty \left[ \sqrt{m+1}(K_{2(m+1)} - K_{2m} ) + (\sqrt{m+1} - \sqrt{m}) K_{2m} \right]
\ee which implies that
\be
\ds \frac{\text{(II)}}{Y_{m+1} - Y_m} = \tau_\infty \left[ \sqrt{m+1} +  (\sqrt{m+1} - \sqrt{m})\ds \frac{1+ \tau_m}{2} \frac{K_{2m}}{ \vp_{2m}^2} \right]
\label{e5}
\ee

Next, we show that $\lim_{m \to \infty} K_{2m}/\vp_{2m}^2 = -1/\tau_\infty$ by the Stolz--Ces\`aro Theorem.
\be \begin{array}{ll}
\ds \lim_{m \to \infty} \frac{K_{2m}}{\vp_{2m}^2} & = \left( \ds \lim_{m \to \infty} \ds \frac{\vp_{2(m+1)}^2 - \vp_{2m}^2}{K_{2(m+1)} - K_{2m}} \right)^{-1} \\
& = \ds \lim_{m \to \infty} \ds \left( {\frac{1-\tau_m}{1+\tau_m} -1}\right)^{-1} \left({\frac{1-\tau_m}{1+\tau_m} +1} \right) \\
&= - \ds \frac 1 {\tau_\infty}
\label{e6}
\end{array}
\ee

%Thus,
%\be
%\frac{(II)}{Y_{m+1} - Y_m} = \tau_\infty \left[ (m+1) -\ds \frac{1+\tau_m}{\tau_\infty + o(1)} \right]
%\label{e6}
%\ee

Combining (\ref{e4}), (\ref{e5}) and (\ref{e6}), we obtain
\be \ds \lim_{m \to \infty} m \left({\Delta_{2m}(1)-(\tau_{\infty})} \right) =1
%= & \ds \lim_{m \to \infty} \frac 1 2 - \tau_\infty\left( m +\frac 1 2 \right) + \ds \frac{m+\frac 1 2 }{m} +\tau_\infty \left[ (m+1) -\ds \frac{1+\tau_m}{2(\tau_\infty + o(1))} \right] \\
%
%\end{array}
\ee

Next, we are going to show that
\be
\Delta_{2m}(1) = -\tau_\infty + \ds \frac 1 {\sqrt m} + o\left( \ds \frac 1 {m}\right)
\label{e13}
\ee by computing the second-order term. We do so by proving that 
\be
L_2 \equiv \ds \lim_{m \to \infty} m \left( \Delta_{2m} -(- \ds \tau_\infty )-\ds \frac 1 {\sqrt m} \right) =0
\ee

Let
\be
P_m = m \vp_{2m}^2 + m \tau_\infty K_{2m} - \sqrt{m} K_{2m}
\ee

Then
\begin{multline}
P_{m+1} - P_m =\left[ (m+1)\ds \frac{1-\tau_m}{1+\tau_m} - m \right] \vp_{2m}^2 \\ + (m+1)\tau_\infty \left[ K_{2(m+1)} - K_{2m}\right]  + \left[ (m+1)-m\right] \tau_\infty K_{2m}\\-\sqrt{m+1} \left[ K_{2(m+1)} - K_{2m}\right] -(\sqrt{m+1}-\sqrt{m})K_{2m}
\end{multline}

Combining with previous results about $Y_{m+1}-Y_m$ and $K_{2m}/\vp_{2m}^2$, we have
\be
L_2 = \ds \lim_{m \to \infty} \frac{P_{m+1} - P_m}{Y_{m+1} - Y_m} = 0
\ee which proves (\ref{e13}).

Next, we will obtain the third-order term by computing
\be
L_3 = \ds \lim_{m \to \infty} m^{3/2} \left( \Delta_{2m} - (-\tau_\infty) - \ds \frac{1}{\sqrt{m}}\right)
\ee

Let
\be
J_m = m^{3/2} \vp_{2m}^2 + m^{3/2} \tau_\infty K_{2m} - m K_{2m}
\ee

By a similar argument as in (\ref{e6}), 
\begin{multline}
J_{m+1} - J_m =\left[ (m+1)^{3/2}\ds \frac{1-\tau_m}{1+\tau_m} - m^{3/2} \right] \vp_{2m}^2 \\ + (m+1)^{3/2} \tau_\infty \left[ K_{2(m+1)} - K_{2m}\right]  + \left[ (m+1)^{3/2}-m^{3/2}\right] \tau_\infty K_{2m}\\-(m+1) \left[ K_{2(m+1)} - K_{2m}\right] -(m+1-m)K_{2m}
\end{multline} which implies that
\be
L_3 = \ds\lim_{m \to \infty} \ds \frac{J_m}{Y_m}=1+\frac{1}{2 \tau_\infty}
\ee

\vspace{2cm}

\subsection{when n is odd} We compute $\Delta_n(1)$ when $n=2m+1$ using the point mass formula (\ref{deltandef}). Let $A_n$ and $B_n$ be defined as in (\ref{eqn61a}) and (\ref{eqn61b}). Then
\be
\lim_{m \to \infty} \Delta_{2m+1}(1) = (1-|\tau_\infty|^2)^{1/2}\lim_{m \to \infty}  \, \ds \frac{A_{2m+1}}{B_{2m+1}}
\ee

We will use the Stolz--Ces\`aro Theorem again. Note that
\be
A_{2(m+1)+1} - A_{2m+1} =\left(  \ds \sqrt{\frac{1-\tau_{m+1}}{1+\tau_{m+1}}} \frac{1-\tau_m}{1+\tau_m}- \ds \sqrt{\frac{1-\tau_{m}}{1+\tau_{m}}} \right) \vp_{2m}^2
\ee and because $\vp_{2m+3}=\vp_{2m+2}$,
\be
B_{2(m+1)+1} - B_{2m+1} = 2 \vp_{2m+2}^2  = 2 \left(\frac{1-\tau_{m}}{1+\tau_{m}} \right) \vp_{2m}^2
\ee

Therefore,
\be
\ds \lim_{m \to \infty} \Delta_{2m+1}(1) = \ds\frac{- \tau_\infty (1-|\tau_\infty|^2)^{1/2} }{\sqrt{(1+\tau_\infty)(1-\tau_\infty)}} = -\tau_\infty
\ee

Next, we prove the rate of convergence by computing
\be
\ds \lim_{m \to \infty} \sqrt m \left({\Delta_{2m+1}(1)+\tau_{\infty}} \right)= 1
\label{e7}
\ee

Since $\alpha_n \in \mathbb{R}$, the recursion relation becomes
\be
(1-|\alpha_n|^2)^{1/2} \vp_{n+1} = \vp_n-\ol{\alpha_n} \vp_n^*= (1-\alpha_n) \vp_n
\ee

Therefore,
\be
\Delta_{2m+1}(1) = \ds \frac{(1-\alpha_{2m+1} )\vp_{2m+1}^2  }{K_{2m+1}} = (1-\tau_{m}) \ds \frac{\vp_{2m}^2}{K_{2m+1}}
\label{e10}
\ee

Let
\begin{eqnarray}
P_m &=& \sqrt{m} \left[ (1-\tau_m) \vp_{2m}^2 + \tau_{\infty} K_{2m+1} \right] \\
Q_m & =& K_{2m+1} \to \infty
\end{eqnarray}

Note that
\be
Q_{m+1}- Q_m = K_{2m+3} - K_{2m+1} = 2\vp_{2(m+1)}^2
\ee and
\begin{multline}
P_{m+1} - P_m = \underbrace{\left[ \sqrt{m+1} (1-\tau_{m+1}) \vp_{2(m+1)}^2 - \sqrt{m}(1-\tau_{m} ) \vp_{2m}^2 \right]}_{\text{(I)}} \\ + \tau_\infty \sqrt{m+1} \left[ K_{2m+3} - K_{2m+1} \right] + \underbrace{(\sqrt{m+1}-\sqrt m) \tau_\infty K_{2m+1}}_{\text{(II)}}
\end{multline}

Since $(1-\tau_m) \vp_{2m}^2 = (1+\tau_m) \vp_{2(m+1)}^2$,
\be \begin{array}{ll}
\ds \frac{\text{(I)}}{Q_{m+1} - Q_m} & = \ds \frac{\sqrt{m+1} ( 1-\tau_{m+1}) - \sqrt m ( 1 + \tau_m)  }{2}
\end{array}
\ee

Next, consider (II). We compute
\be
\ds \lim_{m \to \infty} \ds \frac{\vp_{2(m+1)}^2 - \vp_{2m}^2}{K_{2m+1} -K_{2m-1}}  = \ds \frac{\left( \frac{1-\tau_m}{1+\tau_m} -1\right) \vp_{2m}^2}{2\vp_{2m}^2} = \ds \frac{-\tau_\infty}{1+\tau_\infty}
\ee which implies
\be
\ds \frac{\text{(II)}}{Q_{m+1} - Q_m} = -(1+\tau_\infty)(\sqrt{m+1} - \sqrt m)
\ee

Therefore,
\be
\ds \lim_{m \to \infty} m \left( \Delta_{2m+1}(1) - (-\tau_\infty) \right) = \ds \lim_{m \to \infty} \frac{P_m}{Q_m} = 1
\ee

Next, we will prove that
\be
\Delta_{2m+1}=-\tau_\infty + \ds \frac 1 {\sqrt m} + o\left( \ds \frac 1 m \right)
\label{e15}
\ee
by showing
\be
L_2 '\equiv \ds \lim_{m \to \infty} m \left( \Delta_{2m+1}+\tau_\infty - \frac{1}{\sqrt m} \right) = 0
\ee

As explained in (\ref{e10}), it suffices to consider
\be
H_m = m (1-\tau_m) \vp_{2m}^2 + m\tau_\infty K_{2m+1} - \sqrt{m} K_{2m+1}
\ee

\begin{multline}
H_{m+1} - H_m = \underbrace{(m+1) (1-\tau_{m+1})\vp_{2(m+1)}^2 -m(1-\tau_m)\vp_{2m}^2}_{\text{(I)}} \\ \underbrace{+ (m+1) \tau_\infty K_{2m+3} - m \tau_\infty K_{2m+1}}_{\text{(II)}} 
\underbrace{-\sqrt{m+1} K_{2m+3} + \sqrt{m} K_{2m+1}}_{\text{(III)}}
\end{multline}

Since $(1-\tau_m)\vp_{2m}^2=(1+\tau_m) \vp_{2(m+1)}^2$, we have
\be
\ds\frac{\text{(I)}}{Q_{m+1} - Q_m} = \ds \frac{ (m+1) (1-\tau_{m+1}) -m(1+\tau_m)}{2}
\label{i}
\ee

\be
\ds \frac{\text{(II)}}{Q_{m+1}-Q_m}=\tau_\infty(m+1) +(m+1-m) \tau_\infty \ds \frac{K_{2m+1}}{2 \vp_{2(m+1)}^2}
\label{ii}
\ee

\be
\ds \frac{\text{(III)}}{Q_{m+1}-Q_m}=\ds \frac{(-\sqrt{m+1}+\sqrt{m})K_{2m+3}}{2 \vp_{2(m+1)}^2}+ (-\sqrt{m})
\label{iii}
\ee

This proves that $L'=0$ and thus (\ref{e15}).

Next, we compute
\be
L_3'= \ds \lim_{m \to \infty} m^{3/2} \left( \Delta_{2m+1} + \tau_\infty - \ds \frac{1}{\sqrt{m}} \right)
\ee

By similar arguments as in (\ref{i}), (\ref{ii}) and (\ref{iii}), we conclude that
\be
L_3' = 1+ \ds \frac{1}{2\tau_\infty}
\ee

This concludes the proof of Lemma \ref{lemma1}
\end{proof}

Finally, we apply the Szeg\H o map to this perturbed measure $d \tilde{\mu}_y$ to form the perturbed measure $d\tilde{\gamma}_y$ on $[-2,2]$, which is defined by
\be
\tilde{\gamma}_y(x) = (1-\gamma) d\gamma_y(x) + \gamma \delta_{x=2}
\ee
%\ee where
%\be
%\tilde{w}_y(x) = \pi^{-1} (4-x^2)^{-1/2} f_y(x)
%\ee is the a.c. part of $\tilde{\gamma}_y$.

%By Theorem 13.1.7 of \cite{simon2}, under the Szeg\H o map,
%\begin{eqnarray}
%a_{n+1}^2 & = & (1-\alpha_{2n-1})(1-\alpha_{2n}^2)(1+\alpha_{2n+1}) \\
%b_{n+1}&=& (1-\alpha_{2n-1})\alpha_{2n} - (1+\alpha_{2n-1}) \alpha_{2n-2}
%\end{eqnarray} with the convention that $\alpha_{-1}=1$.

For the sake of convenience, we temporarily denote $\alpha_n \equiv \alpha_n(d\tilde{\mu}_y)$. Since $b_{n}(d\gamma_y)\equiv 0$, if suffices to consider
\be
b_{n+1}(d\gamma_y) =  \ds \frac{1}{\sqrt{n}} - \ds \frac{1}{\sqrt{n-1}} + o\left(\ds \frac 1 {n^{3/2}} \right) = \frac {-1} {2 n^{3/2}} + o\left(\frac 1 {n^{3/2}}\right)
\ee

It is more complicated with $a_n(d\tilde{\gamma_y})$. Recall that
\be
a_{n+1}(d\tilde{\gamma}_y)^2 = (1-\alpha_{2n-1})(1-\alpha_{2n}^2)(1+\alpha_{2n+1}) \\
\ee and we know that
\be
a_{n+1}(d\gamma_y)^2=(1-\tau_{n-1})(1+\tau_n)
\ee

Therefore, upon solving the algebra, we obtain
\be
a_{n+1}(d\tilde{\gamma}_y)^2 - a_{n+1}(d\gamma_y)^2 = \ds \frac{1}{2(1+\tau_\infty) m^{3/2}} + o\left( \ds \frac{1}{m^{3/2}} \right)
\ee

Upon scaling, we have
\be
a_{n+1}^2(d{\gamma}) - a_{n+1}^2(d\gamma_0) =  \ds \frac{2}{y^2 (1+\tau_\infty) m^{3/2}} + o\left( \ds \frac{1}{m^{3/2}} \right)
\ee
\medskip

\section{acknowledgement}
I would like to thank Professor Barry Simon for suggesting this problem and for all the very helpful discussions.

\end{document}